\input amstex
\documentstyle{amsppt}
%
\catcode`@=11
\redefine\output@{%
  \def\break{\penalty-\@M}\let\par\endgraf
  \ifodd\pageno\global\hoffset=105pt\else\global\hoffset=8pt\fi  
  \shipout\vbox{%
    \ifplain@
      \let\makeheadline\relax \let\makefootline\relax
    \else
      \iffirstpage@ \global\firstpage@false
        \let\rightheadline\frheadline
        \let\leftheadline\flheadline
      \else
        \ifrunheads@ 
        \else \let\makeheadline\relax
        \fi
      \fi
    \fi
    \makeheadline \pagebody \makefootline}%
  \advancepageno \ifnum\outputpenalty>-\@MM\else\dosupereject\fi
}
\catcode`\@=\active
\nopagenumbers
\def\negskp{\hskip -2pt}
\def\sign{\operatorname{sign}}
\def\tr{\operatorname{tr}}

\def\Img{\operatorname{Im}}
\def\Ker{\operatorname{Ker}} 
\def\rank{\operatorname{rank}}

\def\blue#1{#1}
\catcode`#=11\def\diez{#}\catcode`#=6
\catcode`_=11\def\podcherkivanie{_}\catcode`_=8
\def\mycite#1{\cite{\blue{#1}}\immediate\special{ps:
     ShrHPSdict begin /ShrBORDERthickness 0 def}}

\def\mytag#1{%
    \tag#1}
\def\mythetag#1{\thetag{\blue{#1}}\immediate\special{ps:
     ShrHPSdict begin /ShrBORDERthickness 0 def}}
\def\myrefno#1{\no#1}
\def\myhref#1#2{\blue{#2}\immediate\special{ps:
     ShrHPSdict begin /ShrBORDERthickness 0 def}}
\def\myEarXivlink{\myhref{http://arXiv.org}{http:/\negskp/arXiv.org}}
\def\myGeoCities{\myhref{http://www.geocities.com}{GeoCities}}
\def\mytheorem#1{\csname proclaim\endcsname{Theorem #1}}
\def\mythetheorem#1{\blue{#1}\immediate\special{ps:
     ShrHPSdict begin /ShrBORDERthickness 0 def}}
\def\mylemma#1{\csname proclaim\endcsname{Lemma #1}}

\def\mycorollary#1{\csname proclaim\endcsname{Corollary #1}}

\pagewidth{360pt}
\pageheight{606pt}
\topmatter
\title
A note on pairs of metrics\\ 
in a three-dimensional linear vector space.
\endtitle
\author
R.~A.~Sharipov
\endauthor
\address 5 Rabochaya street, 450003 Ufa, Russia\newline
\vphantom{a}\kern 12pt Cell Phone: +7(917)476 93 48
\endaddress
\email \vtop to 30pt{\hsize=280pt\noindent
\myhref{mailto:r-sharipov\@mail.ru}
{r-sharipov\@mail.ru}\newline
\myhref{mailto:R\podcherkivanie Sharipov\@ic.bashedu.ru}
{R\_\hskip 1pt Sharipov\@ic.bashedu.ru}\vss}
\endemail
\urladdr
\vtop to 20pt{\hsize=280pt\noindent
\myhref{http://www.geocities.com/r-sharipov}
{http:/\negskp/www.geocities.com/r-sharipov}\newline
\myhref{http://www.freetextbooks.boom.ru/index.html}
{http:/\negskp/www.freetextbooks.boom.ru/index.html}\vss}
\endurladdr
\abstract
    Pairs of metrics in a three-dimensional linear vector 
space are considered, one of which is a Minkowski type
metric with the signature $(+,-,-)$. Such metric pairs are 
classified and canonical presentations for them in each 
class are suggested.
\endabstract
\endtopmatter
\TagsOnRight
\document

\rightheadtext{A note on pairs of metrics \dots}
\head
1. Introduction.
\endhead
     In this paper I continue the study of metric pairs initiated 
in paper \mycite{1}. Let $V$ be a three-dimensional linear vector 
space over the field of real numbers $\Bbb R$. Like in \mycite{1} 
we assume that $V$ is equipped with two metrics $\bold g$ and 
$\check\bold g$. The first metric $\bold g$ is assumed to be a
Minkowski type metric, i\.\,e\. a metric with the signature 
$(+,-,-)$. The actual Minkowski metric arises in the 
four-dimensional space-time of the Special Relativity, while  
$\bold g$ is a three-dimensional model of this Minkowski metric.
\head
2. Associated linear operators. 
\endhead
     The Minkowski type metric $\bold g$ with the signature
$(+,-,-)$ is non-degenerate. For this reason each other metric 
$\check\bold g$ in $V$ produces a linear operator $\check\bold F$
associated with it through the metric $\bold g$. The operator
$\check\bold F$ is defined by the formula
$$
\hskip -2em
g(\check\bold F(\bold X),\bold Y)=\check g(\bold X,\bold Y).
\mytag{2.1}
$$
Due to the symmetry of the form $\check g(\bold X,\bold Y)$ this
formula \mythetag{2.1} is extended to 
$$
\hskip -2em
g(\check\bold F(\bold X),\bold Y)=\check g(\bold X,\bold Y)=
g(\bold X,\check\bold F(\bold Y)).
\mytag{2.2}
$$
The formula \mythetag{2.2} means that the associated operator 
$\check\bold F$ is symmetric with respect to the metric $\bold g$.
It is easy to show that $\check\bold F$ is symmetric with respect
to $\check\bold g$ as well:
$$
\check g(\check\bold F(\bold X),\bold Y)
=g(\check\bold F(\bold X),\check\bold F(\bold Y))
=\check g(\bold X,\check\bold F(\bold Y)).
$$\par
     If the metrics $\bold g$ and $\check\bold g$ are given by their
components $g_{ij}$ and $\check g_{ij}$ in some basis $\bold e_0,\,
\bold e_1,\,\bold e_2$, then the associated operator $\check\bold F$ 
is given by the components
$$
\hskip -2em
\check F^i_j=\sum^2_{s=0} g^{is}\,\check g_{sj}.
\mytag{2.3}
$$
Here $g^{ij}$ are the components of the metric dual for $\bold g$.
They form the matrix inverse to the matrix formed by the components 
$g_{ij}$. Let $P(\lambda)$ be the characteristic polynomial of 
the operator $\check\bold F$. Then we have
$$
\hskip -2em
P(\lambda)=\det(\check\bold F-\lambda\,\bold I)
=-\lambda^3+a_0\,\lambda^2-a_1\,\lambda+a_2.
\mytag{2.4}
$$
The coefficients $a_0$, $a_1$, $a_2$ of the characteristic polynomial 
\mythetag{2.4} are the invariants of the pair of metrics $\bold g$ and
$\check\bold g$. They are calculated as follows:
$$
\xalignat 3
&\hskip -2em
a_0=\tr\check\bold F,
&&a_1=\frac{(\tr\check\bold F)^2-\tr\check\bold F^2}{2},
&&a_2=\det\check\bold F.
\qquad
\mytag{2.5}
\endxalignat
$$\par
     A third order polynomial with real coefficients has at least one
root in the field of real numbers $\Bbb R$. This real root is an eigenvalue 
of the associated operator $\check\bold F$. It corresponds to at least one 
eigenvector $\bold v\in V$. For this reason we consider the following three 
mutually exclusive cases:
\roster
\rosteritemwd=2pt 
\item there is a simple eigenvalue $\lambda=\lambda_0$ of $\check\bold F$ 
with a time-like eigenvector $\bold v$ with respect to the metric $\bold g$,
 i\.\,e\. such that $g(\bold v,\bold v)>0$;
\item there are no simple eigenvalues with time-like eigenvectors, but 
there is a simple eigenvalue $\lambda=\lambda_2$ of $\check\bold F$ with 
a space-like eigenvector $\bold v$ with respect to the metric $\bold g$, 
i\.\,e\. such that $g(\bold v,\bold v)<0$;
\item there are no simple eigenvalues of $\check\bold F$ at all.
\endroster
     Saying a {\bf simple eigenvalue} I mean a real root of the polynomial
\mythetag{2.4} whose multiplicity $k=1$. In the third case, where there no 
simple eigenvalues, a real root of the polynomial \mythetag{2.4}, which does 
always exist, is unique and it is a triple eigenvalue with the multiplicity 
$k=3$. The subdivision of our consideration into the above three cases is 
based on the following theorem.
\mytheorem{2.1} The associated operator $\check\bold F$ of a metric pair 
$\bold g$ and $\check\bold g$ such that $\bold g$ is a Minkowski type metric 
with the signature $(+,-,-)$ cannot have a simple eigenvalue with an
eigenvector $\bold v$ such that $g(\bold v,\bold v)=0$.
\endproclaim
\demo{Proof} For to prove the theorem~\mythetheorem{2.1} we consider
the complexification $\Bbb CV=\Bbb C\otimes V$ of the vector space 
$V$. It is naturally equipped with the complexifications of the metrics 
$\bold g$ and $\check\bold g$. Their associated operator $\check \bold F$ 
coincides with the complexification of the operator $\check\bold F$ acting 
in $V$. The complexified metrics in $\Bbb CV$ inherit their signatures 
from $\bold g$ and $\check\bold g$ with the only difference --- all minuses
turn to pluses. In particular, the signature of the complexified metric
$\bold g$ is $(+,+,+)$.\par
     The complex linear vector space $\Bbb CV=\Bbb C\otimes V$ is 
naturally equipped with the semilinear involution of complex conjugation:
$$
\hskip -2em
\tau\!:\,\Bbb CV\to\Bbb CV.
\mytag{2.6}
$$
The space $V$ is embedded into $\Bbb CV$ as an $\Bbb R$-linear subspace 
invariant under the involution \mythetag{2.6}. Assume that the associated 
operator $\check\bold F$ has a simple eigenvalue $\lambda_0\in\Bbb R$ 
with an eigenvector $\bold v_0\in V\subset\Bbb CV$ such that 
$$
\hskip -2em
g(\bold v_0,\bold v_0)=0.
\mytag{2.7}
$$
Apart from $\lambda_0$, the operator $\check\bold F$ has two distinct
eigenvalues $\lambda_1\neq\lambda_2$ or one double eigenvalue $\lambda_1
=\lambda_2$. If $\lambda_1\neq\lambda_2$ we have two extra eigenvectors
$\bold v_1$ and $\bold v_2$ in $\Bbb CV$. The eigenvectors $\bold v_0$,
$\bold v_1$, and $\bold v_2$, corresponding to three distinct eigenvalues
$\lambda_0$, $\lambda_1$ and $\lambda_2$ are linearly independent (see
\mycite{2}). They form a basis in $\Bbb CV$. Moreover, the eigenvectors
$\bold v_i$ and $\bold v_j$ of the associated operator $\check\bold F$
corresponding to $\lambda_i\neq\lambda_j$ are orthogonal with respect to
the metric $\bold g$. Indeed, from \mythetag{2.2} we derive
$$
\hskip -2em
\lambda_i\,g(\bold v_i,\bold v_j)
=g(\check\bold F(\bold v_i),\bold v_j)
=g(\bold v_i,\check\bold F(\bold v_j))
=\lambda_j\,g(\bold v_i,\bold v_j).
\mytag{2.8}
$$
Since $\lambda_i\neq\lambda_j$, the equality \mythetag{2.8} yields
$$
\hskip -2em
g(\bold v_i,\bold v_j)=0.
\mytag{2.9}
$$
The formula \mythetag{2.9} means that the metric $\bold g$ is 
diagonalized in the basis formed by the eigenvectors $\bold v_0$,
$\bold v_1$, $\bold v_2$, while the formula \mythetag{2.7} says
that it is a degenerate metric with at least one zero in its
signature. This result contradicts the initial assumption that
the extension of $\bold g$ to $\Bbb CV$ is a metric with the
signature $(+,+,+)$.\par
     If $\lambda_1=\lambda_2$ is a double eigenvalue, we consider 
the operator $\bold h=(\check\bold F-\lambda_1\,\bold I)^2$. The
kernel of this operator $W=\Ker\bold h$ is the two-dimensional root
subspace (see \mycite{2}) corresponding to the double eigenvalue 
$\lambda_1=\lambda_2$. In this case we have
$$
\hskip -2em
\Bbb CV=\bigl<\bold v_0\bigr>\oplus W,
\mytag{2.10}
$$
where $\bigl<\bold v_0\bigr>$ is the linear span of the eigenvector 
$\bold v_0$. The subspaces $\bigl<\bold v_0\bigr>$ and $W$ in 
\mythetag{2.10} are perpendicular to each other. Indeed, let
$\bold w\in W$. Then
$$
(\lambda_0-\lambda_1)^2\,g(\bold v_0,\bold w)
=g(\bold h(\bold v_0),\bold w)
=g(\bold v_0,\bold h(\bold w))=0
\mytag{2.11}
$$
since $\bold h(\bold w)=0$ by the definition of a root subspace
(see \mycite{2}). From \mythetag{2.11}, using $\lambda_0\neq\lambda_1$,
we derive $\bold v_0\perp W$, while from $\bold v_0\perp W$ we
conclude that $\bold g$ should have at least one zero in its 
signature. Again, this result is a contradiction to the initial
assumption that the complexification of $\bold g$ is a metric with 
the signature $(+,+,+)$ in $\Bbb CV$. The theorem~\mythetheorem{2.1} 
is proved.\qed\enddemo
\head
3. The first case.
\endhead
     In the previous section we have divided the study of metric pairs 
into three mutually exclusive cases. In the {\bf first case} the
associated operator $\check\bold F$ has an eigenvalue $\lambda_0$
whose eigenvector $\bold v_0$ is a time-like vector with respect
to the metric $\bold g$, i\.\,e. $g(\bold v_0,\bold v_0)>0$. We can 
normalize this eigenvector so that
$$
\hskip -2em
g(\bold v_0,\bold v_0)=1.
\mytag{3.1}
$$  
Let's denote by $W$ the orthogonal complement of $\bold v_0$ with
respect to the metric $\bold g$:
$$
\hskip -2em
W=\{\bold w\in V\!:\ g(\bold v_0,\bold w)=0\}.
\mytag{3.2}
$$
Since $\bold g$ is non-degenerate, \mythetag{3.2} is a two-dimensional 
subspace in $V$. From \mythetag{3.1} we derive that $W$ is transversal
to the vector $\bold v_0$:
$$
V=\bigl<\bold v_0\bigr>\oplus W.
\mytag{3.3}
$$
The expansion \mythetag{3.3} is similar to the expansion \mythetag{2.10}. 
Using \mythetag{2.2} and \mythetag{3.2}, we easily prove that the 
subspace $W$ is perpendicular to $\bold v_0$ with respect to the metric 
$\check\bold g$ as well. Indeed, if $\bold w\in W$, then we have
$$
\hskip -2em
\check g(\bold v_0,\bold w)
=g(\check\bold F(\bold v_0),\bold w)
=\lambda_0\,g(\bold v_0,\bold w)=0.
\mytag{3.4}
$$
Thus the expansion \mythetag{3.3} is an expansion of $V$ into a
direct sum of two subspaces mutually perpendicular with respect to
both metrics $\bold g$ and $\check\bold g$.\par
Due to \mythetag{3.1}
the restriction of $\bold g$ to $W$ is purely negative. Therefore 
the restrictions of $\bold g$ and $\check\bold g$ to $W$ can be
diagonalized simultaneously. If the vectors $\bold v_1$ and 
$\bold v_2$ form a basis of $W$ where both of these restrictions are 
diagonal, then we can complement them with the vector $\bold v_0$. 
As a result we get a basis $\bold v_0,\bold v_1,\,\bold v_2$ in
$V$ such that the metrics $\bold g$ and $\check\bold g$ are given
by the following matrices in this basis:
$$
\xalignat 2
&\hskip -2em
g_{ij}=\Vmatrix\format\l\ \ &\r &\ \ \r\\
1 & 0 & 0\\ 
\vspace{1ex} 0 &-1 & 0\\ 
\vspace{1ex} 0 & 0 &-1\endVmatrix,
&&\check g_{ij}=\Vmatrix\format\l\ \ \ &\r &\ \ \ \r\\
a & 0 & 0\\ 
\vspace{1ex} 0 &b & 0\\ 
\vspace{1ex} 0 & 0 &c\endVmatrix.
\mytag{3.5}
\endxalignat
$$
Applying the formula \mythetag{2.3} to the components of the matrices
\mythetag{3.5}, we derive that the associated operator $\check\bold F$ 
is given by the matrix
$$
\hskip -2em
F^i_j=\Vmatrix\format\l\ \ \ &\r &\ \ \ \r\\
a & 0 & 0\\ 
\vspace{1ex} 0 &-b & 0\\ 
\vspace{1ex} 0 & 0 &-c\endVmatrix
\mytag{3.6}
$$
in the basis $\bold v_0,\bold v_1,\,\bold v_2$. From \mythetag{3.6}
we find that $\lambda_0=a$. Therefore, we have 
$$
\xalignat 2
&\hskip -2em
a\neq -b,
&&a\neq -c.
\mytag{3.7}
\endxalignat
$$\par
\parshape 18 170pt 190pt 170pt 190pt 170pt 190pt 170pt 190pt
 170pt 190pt 170pt 190pt 170pt 190pt 170pt 190pt 170pt 190pt 
 170pt 190pt 170pt 190pt 170pt 190pt 170pt 190pt 170pt 190pt 
 170pt 190pt 170pt 190pt 170pt 190pt 0pt 360pt
     Using the matrix \mythetag{3.6} and taking into account 
\mythetag{3.7}, we can calculate the characteristic polynomial 
\mythetag{2.4} explicitly:
$$
P(\lambda)=-(\lambda-a)(\lambda+b)(\lambda+c).
\qquad
\mytag{3.8}
$$
In general case, apart from \mythetag{3.7}, we have the following
inequality:
$$
\hskip -2em
b\neq c.
\mytag{3.9}
$$
In this general case the graph of the polynomial \mythetag{3.8} 
intersects the $\lambda$-axis at three distinct points as shown 
on Fig\.~3.1. \vadjust{\vskip 5pt\hbox to 0pt{\kern -10pt
\includegraphics{ThreeMetric01.eps}\hss}\vskip -5pt}Now
we shall derive the condition providing the polynomial 
\mythetag{2.4} to have three simple roots. Let's denote by
$Q(\lambda)$ \pagebreak the derivative of the polynomial 
\mythetag{2.4}: $Q(\lambda)=P'(\lambda)$. Here is the explicit 
formula for this derivative:
$$
\hskip -2em
Q(\lambda)=-3\,\lambda^2+2\,a_0\,\lambda-a_1.
\mytag{3.10}
$$
Let's denote by $D_2$ the discriminant of the polynomial 
\mythetag{3.10}. Then we have
$$
\hskip -2em
D_2=4\,(a_0)^2-12\,a_1.
\mytag{3.11}
$$
Using the discriminant \mythetag{3.11}, we calculate $\lambda_{\text{max}}$
and $\lambda_{\text{min}}$ on Fig\.~3.1:
$$
\xalignat 2
&\hskip -2em
\lambda_{\text{min}}=\frac{a_0}{3}-\frac{\sqrt{\mathstrut D_2}}{6},
&&\lambda_{\text{max}}=\frac{a_0}{3}+\frac{\sqrt{\mathstrut D_2}}{6}.
\mytag{3.12}
\endxalignat
$$
It is clear that the characteristic polynomial \mythetag{2.4} has 
three distinct real roots if and only if the following conditions 
are fulfilled:
$$
\xalignat 2
&\hskip -2em
P(\lambda_{\text{min}})<0,
&&P(\lambda_{\text{max}})>0.
\mytag{3.13}
\endxalignat
$$
Substituting \mythetag{3.12} into \mythetag{2.4}, by means of direct 
calculations we find that the inequalities \mythetag{3.13} are equivalent 
to the following ones:
$$
\hskip -2em
\aligned
&(D_2)^{3/2}>8\,(a_0)^3-36\,a_0\,a_1+108\,a_2,\\
&(D_2)^{3/2}>-8\,(a_0)^3+36\,a_0\,a_1-108\,a_2.
\endaligned
\mytag{3.14}
$$
The inequalities \mythetag{3.14} can be united into one inequality: 
$$
\hskip -2em
(D_2)^{3/2}>|8\,(a_0)^3-36\,a_0\,a_1+108\,a_2|.
\mytag{3.15}
$$
Squaring both sides of \mythetag{3.15}, we get
$$
\hskip -2em
(D_2)^3>\left(8\,(a_0)^3-36\,a_0\,a_1+108\,a_2\right)^2.
\mytag{3.16}
$$
In particular, the inequality \mythetag{3.16} yields $D_2>0$. For this reason
we did not write $D_2>0$ as a separate condition along with \mythetag{3.13}.
\par
     Let's substitute \mythetag{3.11} into \mythetag{3.16}. As a result, the
above inequality \mythetag{3.16} is transformed to the following inequality:
$$
\hskip -2em
-27\,(a_2)^2+18\,a_0\,a_1\,a_2+(a_1)^2\,(a_0)^2-4\,(a_0)^3\,a_2-4\,(a_1)^3>0.
\mytag{3.17}
$$
As appears, the left hand side of the inequality \mythetag{3.17} coincides
with the discriminant of the polynomial \mythetag{2.4} itself, i\.\,e\. we
have
$$
D_3=-27\,(a_2)^2+18\,a_0\,a_1\,a_2+(a_1)^2\,(a_0)^2-4\,(a_0)^3\,a_2-4\,(a_1)^3.
\quad
\mytag{3.18}
$$
Now the inequality \mythetag{3.17} is written as
$$
\hskip -2em
D_3>0.
\mytag{3.19}
$$
\mytheorem{3.1}The cubic polynomial \mythetag{2.4} with the real 
coefficients $a_0$, $a_1$, $a_2$ has three distinct real roots if 
and only if its discriminant \mythetag{3.18} is positive.
\endproclaim
     The inequality \mythetag{3.19} is the only condition that distinguishes
the case of three simple roots within the first subcase in the classification 
scheme suggested in section~2. Indeed, if the associated operator $\check
\bold F$ has three simple eigenvalues, then we have three linearly independent
eigenvectors mutually perpendicular with respect to the metric $\bold g$.
One of them is certainly a time-like vector, while two others are space-like
vectors according to the signature $(+,-,-)$ of the metric $\bold g$.\par
     Apart from the general case, there is a special subcase within the 
first case. In this special subcase the inequality \mythetag{3.9} is broken 
and we have the equality
$$
\hskip -2em
b=c.
\mytag{3.20}
$$
Due to \mythetag{3.20} the graph of the characteristic polynomial 
\mythetag{2.4} takes one of two possible shapes shown on Fig\.~3.2
and on Fig\.~3.3. \vadjust{\vskip 5pt\hbox to 0pt{\kern 5pt
\includegraphics{ThreeMetric02.eps}\hss}\vskip 205pt}Due 
to \mythetag{3.20} the inequality \mythetag{3.19} for the discriminant
$D_3$ in this special subcase turns to the equality
$$
\hskip -2em
D_3=0.
\mytag{3.21}
$$
The inequality \mythetag{3.15} is replaced by the equality 
$$
\hskip -2em
(D_2)^3=\left(8\,(a_0)^3-36\,a_0\,a_1+108\,a_2\right)^2,
\mytag{3.22}
$$
which is equivalent to the equality \mythetag{3.21}. The equality 
\mythetag{3.22} should be complemented with the inequality for the 
discriminant $D_2$:
$$
\hskip -2em
D_2>0.
\mytag{3.23}
$$
The inequality \mythetag{3.23} follows from the inequalities 
\mythetag{3.7}.\par
     Note that the equality \mythetag{3.22} leads to one of the
following two equalities resembling the above inequalities 
\mythetag{3.14}:
$$
\hskip -2em
\aligned
&(D_2)^{3/2}=8\,(a_0)^3-36\,a_0\,a_1+108\,a_2,\\
&(D_2)^{3/2}=-8\,(a_0)^3+36\,a_0\,a_1-108\,a_2.
\endaligned
\mytag{3.24}
$$
The equalities \mythetag{3.24} are mutually exclusive. If the
first equality \mythetag{3.24} holds, we have $P(\lambda_{\text{min}})
=0$ and the graph of the polynomial $P(\lambda)$ takes the shape 
presented on Fig\.~3.2. Otherwise, if the second equality \mythetag{3.24} 
holds, then $P(\lambda_{\text{max}})=0$ and the graph of the polynomial 
$P(\lambda)$ takes the shape presented on Fig\.~3.3. Applying
\mythetag{3.11} to \mythetag{3.24}, we obtain the formulas
$$
\hskip -2em
\aligned
&\sqrt{\mathstrut D_2}=\frac{2\,(a_0)^3-9\,a_0\,a_1+27\,a_2}
{(a_0)^2-3\,a_1},\\
\vspace{1.5ex}
&\sqrt{\mathstrut D_2}=-\frac{2\,(a_0)^3-9\,a_0\,a_1+27\,a_2}
{(a_0)^2-3\,a_1}
\endaligned
\mytag{3.25}
$$
for these two cases. Substituting \mythetag{3.25} into \mythetag{3.12}, 
we derive the formula
$$
\hskip -2em
\lambda_1=\frac{a_0}{3}-\frac{2\,(a_0)^3-9\,a_0\,a_1+27\,a_2}
{6\,(a_0)^2-18\,a_1}.
\mytag{3.26}
$$
Here $\lambda_1=\lambda_{\text{min}}$ or $\lambda_1=\lambda_{\text{max}}$
depending on which equality holds in \mythetag{3.24}. The formula 
\mythetag{3.26} is a formula for the double root of the polynomial
$P(\lambda)$. It is valid provided the conditions \mythetag{3.21} and
\mythetag{3.23} are fulfilled. Note that
$$
\hskip -2em
a_0=\tr\check\bold F=\lambda_0+\lambda_1+\lambda_2
=\lambda_0+2\,\lambda_1
\mytag{3.27}
$$
since $\lambda_1=\lambda_2$. From \mythetag{3.26} and \mythetag{3.27}
we derive
$$
\hskip -2em
\lambda_0=\frac{a_0}{3}+\frac{2\,(a_0)^3-9\,a_0\,a_1+27\,a_2}
{3\,(a_0)^2-9\,a_1}.
\mytag{3.28}
$$
The formula \mythetag{3.28} is a formula for the simple root of the 
polynomial $P(\lambda)$. It is also valid provided the conditions 
\mythetag{3.21} and \mythetag{3.23} are fulfilled.\par
     Thus, if $\lambda_1=\lambda_2\neq\lambda_0$, the simple eigenvalue
$\lambda_0$ of the associated operator $\check\bold F$ is expressed 
through the invariants \mythetag{2.5} according to the formula 
\mythetag{3.28}. Let $\bold v_0$ be an eigenvector of the operator
$\check\bold F$ corresponding to the eigenvalue $\lambda_0$. Let's
denote
$$
\hskip -2em
\sigma_0=\sign(g(\bold v_0,\bold v_0)).
\mytag{3.29}
$$
The quantity \mythetag{3.29} is a special invariant of a pair of 
metrics $\bold g$ and $\check\bold g$. According to the 
theorem~\mythetheorem{2.1} the invariant $\sigma_0$ can take only
two possible values $\sigma_0=-1$ and $\sigma_0=1$. In the special 
subcase of the first case we have
$$
\hskip -2em
\sigma_0=1.
\mytag{3.30}
$$
This special subcase in our classification scheme is completely described 
by the inequality \mythetag{3.23} and by the equalities \mythetag{3.21} 
and \mythetag{3.30}.
\head
4. The second case.
\endhead
     According to our classification scheme in section~2, in the
{\bf second case} the associated operator $\check\bold F$ has no 
simple eigenvalues with time-like eigenvectors, but it has a simple 
eigenvalue with a space-like eigenvector. \pagebreak Let's denote 
this eigenvalue by $\lambda_2$ and its space-like eigenvector by 
$\bold v_2$. We can normalize $\bold v_2$ by the condition
$$
\hskip -2em
g(\bold v_2,\bold v_2)=-1.
\mytag{4.1}
$$  
Let's denote by $W$ the orthogonal complement to $\bold v_2$ with
respect to the metric $\bold g$:
$$
\hskip -2em
W=\{\bold w\in V\!:\ g(\bold v_2,\bold w)=0\}.
\mytag{4.2}
$$
It is clear that $W$ is a two-dimensional subspace in $V$ transversal
to the vector $\bold v_2$:
$$
V=W\oplus \bigl<\bold v_2\bigr>.
\mytag{4.3}
$$
The subspace $W$ in \mythetag{4.3} is perpendicular to the vector
with respect to both metrics $\bold g$ and $\check\bold g$. The
arguments here are the same as in the case of \mythetag{3.4}.\par
     Let's consider the restrictions of the metrics $\bold g$ and
$\check\bold g$ to the subspace \mythetag{4.2}. The restriction 
of $\bold g$ to $W$ is a metric with the signature $(+,-)$. For 
this reason we can apply the results of \mycite{1} to the pair
of restricted metrics in $W$. The classification scheme of
\mycite{1} includes five cases. Only four of them are applicable 
here. These four cases are those where the restricted metrics 
$\bold g$ and $\check\bold g$ cannot be simultaneously diagonalized
in $W$.\par
     The {\bf subcase one} is fixed by the condition that the
associated operator of the restricted metrics $\bold g$ and
$\check\bold g$ in $W$ has no real eigenvalues. In terms of the 
invariants of the associated operator $\check\bold F$ in $V$ 
this condition is written as the inequality 
$$
\hskip -2em
D_3<0.
\mytag{4.4}
$$
\vskip -1ex
\mytheorem{4.1} The cubic polynomial \mythetag{2.4} with the real 
coefficients $a_0$, $a_1$, $a_2$ has three distinct roots one of which
is real and two others are complex numbers if and only if its discriminant 
\mythetag{3.18} is negative.
\endproclaim
\noindent Under the condition \mythetag{4.4} the matrices of the metrics
$\bold g$ and $\check\bold g$ are brought to 
$$
\xalignat 2
&\hskip -2em
g_{ij}=\Vmatrix\format\l\ \ &\r &\ \ \r\\
1 & 0 & 0\\ 
\vspace{1ex} 0 &-1 & 0\\ 
\vspace{1ex} 0 & 0 &-1\endVmatrix,
&&\check g_{ij}=\Vmatrix\format\l\ \ &\r &\,\ \ \ \r\\
a & b & 0\\ 
\vspace{1ex} b &-a & 0\\ 
\vspace{1ex} 0 & 0 &c\endVmatrix.
\mytag{4.5}
\endxalignat
$$
\mytheorem{4.2} If \,the inequality \mythetag{4.4} is fulfilled, then
the matrices of the metrics $\bold g$ and $\check\bold g$ take their
canonical forms \mythetag{4.5} with $b\neq 0$ in some basis.
\endproclaim
     The {\bf subcase two} is fixed by the condition that the
associated operator of the restricted metrics $\bold g$ and
$\check\bold g$ in $W$ has a double eigenvalue, which is a real 
number, and the $\sigma$-invariant of the restricted metrics
is equal to zero (see \mycite{1}). Let's denote by $\sigma_1$
the $\sigma$-invariant of the restricted metrics in $W$. Then 
we have
$$
\xalignat 4
\hskip -2em
&D_3=0,
&&D_2\neq 0,
&&\sigma_0=-1,
&&\sigma_1=0.
\quad
\mytag{4.6}
\endxalignat
$$
The condition $\sigma_0=-1$ in \mythetag{4.6} is derived from 
\mythetag{4.1}. Under the conditions \mythetag{4.6} the matrices 
of the metrics $\bold g$ and $\check\bold g$ are brought to 
$$
\xalignat 2
&\hskip -2em
g_{ij}=\Vmatrix\format\l\ \ &\r &\ \ \r\\
1 & 0 & 0\\ 
\vspace{1ex} 0 &-1 & 0\\ 
\vspace{1ex} 0 & 0 &-1\endVmatrix,
&&\check g_{ij}=\Vmatrix\format\l\ \ &\r &\,\ \ \ \r\\
a & 0 & 0\\ 
\vspace{1ex} 0 &-a & 0\\ 
\vspace{1ex} 0 & 0 &c\endVmatrix.
\mytag{4.7}
\endxalignat
$$
\mytheorem{4.3} If \,the conditions \mythetag{4.6} are fulfilled, then
the matrices of the metrics $\bold g$ and $\check\bold g$ take their
canonical forms \mythetag{4.7} with $a\neq -c$ in some basis.
\endproclaim
     The {\bf subcase three} differs from the subcase two only by 
the value of the invariant $\sigma_1$. It is fixed by the following 
conditions:
$$
\xalignat 4
\hskip -2em
&D_3=0,
&&D_2\neq 0,
&&\sigma_0=-1,
&&\sigma_1=1.
\quad
\mytag{4.8}
\endxalignat
$$
Under the condition \mythetag{4.8} the matrices of the metrics
$\bold g$ and $\check\bold g$ are brought to 
$$
\xalignat 2
&\hskip -2em
g_{ij}=\Vmatrix\format\l\ \ \ \ &\r &\ \ \r\\
0 & 1 & 0\\ 
\vspace{1ex} 1 &0 & 0\\ 
\vspace{1ex} 0 & 0 &-1\endVmatrix,
&&\check g_{ij}=\Vmatrix\format\l\ \ \ \ &\r &\,\ \ \ \r\\
1 & a & 0\\ 
\vspace{1ex} a &0 & 0\\ 
\vspace{1ex} 0 & 0 &c\endVmatrix.
\mytag{4.9}
\endxalignat
$$
\mytheorem{4.4} If \,the conditions \mythetag{4.8} are fulfilled, then
the matrices of the metrics $\bold g$ and $\check\bold g$ take their
canonical forms \mythetag{4.9} with $a\neq -c$ in some basis.
\endproclaim
     The {\bf subcase four} also differs from the subcase two only by 
the value of the invariant $\sigma_1$. It is fixed by the following 
conditions:
$$
\xalignat 4
\hskip -2em
&D_3=0,
&&D_2\neq 0,
&&\sigma_0=-1,
&&\sigma_1=-1.
\quad
\mytag{4.10}
\endxalignat
$$
Under the condition \mythetag{4.10} the matrices of the metrics
$\bold g$ and $\check\bold g$ are brought to 
$$
\xalignat 2
&\hskip -2em
g_{ij}=\Vmatrix\format\l\ \ \ \ &\r &\ \ \r\\
0 & 1 & 0\\ 
\vspace{1ex} 1 &0 & 0\\ 
\vspace{1ex} 0 & 0 &-1\endVmatrix,
&&\check g_{ij}=\Vmatrix\format\l\ \ \,&\r &\,\ \ \ \r\\
0 & a & 0\\ 
\vspace{1ex} a &-1 & 0\\ 
\vspace{1ex} 0 & 0 &c\endVmatrix.
\mytag{4.11}
\endxalignat
$$
\mytheorem{4.5} If \,the conditions \mythetag{4.10} are fulfilled, then
the matrices of the metrics $\bold g$ and $\check\bold g$ take their
canonical forms \mythetag{4.11} with $a\neq -c$ in some basis.
\endproclaim
\head
5. The third case.
\endhead
    In the {\bf third case} of our classification scheme from the
section~2 the associated operator $\check\bold F$ has exactly one 
real eigenvalue	$\lambda_0$ of the multiplicity $k=3$. This case
is specified by the following equalities:
$$
\xalignat 2
&\hskip -2em
D_3=0,
&&D_2=0.
\mytag{5.1}
\endxalignat
$$
The unique eigenvalue $\lambda_0$ in the third case is expressed 
through the invariant $a_0$:
$$
\hskip -2em
\lambda_0=\frac{a_0}{3}.
\mytag{5.2}
$$
The formula \mythetag{5.2} is analogous to the formulas \mythetag{3.26}
and \mythetag{3.28}. In the third case we define the following invariant 
of the associated operator:
$$
\hskip -2em
\sigma_2=\rank(\check\bold F-\lambda_0\,\bold I).
\mytag{5.3}
$$\par
     The {\bf subcase one} within the third case is determined by the 
following condition for the value of the integer invariant $\sigma_2$ 
in \mythetag{5.3}:
$$
\hskip -2em
\sigma_2=0.
\mytag{5.4}
$$
Combining the conditions \mythetag{5.1} with the condition \mythetag{5.4}, 
we get 
$$
\xalignat 3
&\hskip -2em
D_3=0,
&&D_2=0,
&&\sigma_2=0.
\mytag{5.5}
\endxalignat
$$
Under the condition \mythetag{5.4} the associated operator $\check\bold F$
is a scalar operator:
$$
\hskip -2em
\check\bold F=\lambda_0\,\bold I.
\mytag{5.6}
$$
Due to \mythetag{5.6}, diagonalizing the metric $\bold g$, we simultaneously
diagonalize the second metric $\check\bold g$. For this reason the matrices
of these metrics are brought to
$$
\xalignat 2
&\hskip -2em
g_{ij}=\Vmatrix\format\l\ \ &\r &\ \ \r\\
1 & 0 & 0\\ 
\vspace{1ex} 0 &-1 & 0\\ 
\vspace{1ex} 0 & 0 &-1\endVmatrix,
&&\check g_{ij}=\Vmatrix\format\l\ \ &\r &\ \ \r\\
a & 0 & 0\\ 
\vspace{1ex} 0 &-a & 0\\ 
\vspace{1ex} 0 & 0 &-a\endVmatrix,
\mytag{5.7}
\endxalignat
$$
where $a=\lambda_0$. This result is formulated as a theorem.
\mytheorem{5.1} If \,the conditions \mythetag{5.5} are fulfilled, then
the matrices of the metrics $\bold g$ and $\check\bold g$ take their
canonical forms \mythetag{5.7} in some basis.
\endproclaim
    Now, keeping the conditions \mythetag{5.1} unchanged, we increase
by one the value of the integer invariant $\sigma_2$ in \mythetag{5.4}.
As a result we get
$$
\hskip -2em
\sigma_2=1.
\mytag{5.8}
$$
If the equality \mythetag{5.8} is fulfilled, we have the following two
subspaces in $V$:
$$
\xalignat 2
&\hskip -2em
W=\Img(\check\bold F-\lambda_0\,I),
&&U=\Ker(\check\bold F-\lambda_0\,I).
\mytag{5.9}
\endxalignat
$$
For the dimensions of the subspaces \mythetag{5.9} we have
$$
\xalignat 2
&\hskip -2em
\dim W=1,&&\dim U=2.
\mytag{5.10}
\endxalignat
$$
Moreover, in addition to \mythetag{5.10} we have the following inclusions:
$$
\hskip -2em
W\subset U\subset V.
\mytag{5.11}
$$
Let's choose some nonzero vector $\bold e_0\in W$. Due to \mythetag{5.10} 
this vector is unique up to some nonzero numeric factor. Since $\bold e_0
\in W$ and $W\subset U$, from \mythetag{5.9} we derive
$$
\check\bold F(\bold e_0)=\lambda_0\,\bold e_0.
\mytag{5.12}
$$
The equality \mythetag{5.12} means that $\bold e_0$ is an eigenvector
of the associated operator $\check\bold F$ corresponding to its
eigenvalue \mythetag{5.2}. By the definition of the subspace $W$ in 
\mythetag{5.9} there is another nonzero vector $\bold e_1$ such that
$$
\hskip -2em
\bold e_0=\check\bold F(\bold e_1)-\lambda_0\,\bold e_1.
\mytag{5.13}
$$
Apart from $\bold e_0$ and $\bold e_1$, we choose some nonzero vector 
$\bold e_2\in U$ such that $\bold e_0,\,\bold e_2$ is a basis of the
two-dimensional subspace $U$. Since $\bold e_2\in U$, it is another
eigenvector of the operator $\check\bold F$. Indeed, from \mythetag{5.9} 
we derive
$$
\check\bold F(\bold e_2)=\lambda_0\,\bold e_2.
\mytag{5.14}
$$
Since $\bold e_0\neq 0$, from \mythetag{5.13} we conclude that $\bold e_1
\notin U$. Then from \mythetag{5.11} we derive that the vectors $\bold e_0$, 
$\bold e_1$, and $\bold e_2$ are linearly independent. They form a basis in 
$V$. Let's study the components of the metric $\bold g$ in this basis. 
Applying the formulas \mythetag{5.12}, \mythetag{5.13}, and 
\mythetag{2.2}, we derive
$$
\gathered
g_{00}=g(\bold e_0,\bold e_0)=g(\bold e_0,\check\bold F(\bold e_1))
-\lambda_0\,g(\bold e_0,\bold e_1)=\\
\vspace{1ex}
=g(\check\bold F(\bold e_0),\bold e_1)-\lambda_0\,g(\bold e_0,\bold e_1)
=\lambda_0\,g(\bold e_0,\bold e_1)-\lambda_0\,g(\bold e_0,\bold e_1)=0.
\endgathered
\mytag{5.15}
$$
Similarly, applying the formulas \mythetag{5.14}, \mythetag{5.13}, and 
\mythetag{2.2}, we derive
$$
\gathered
g_{02}=g_{20}=g(\bold e_2,\bold e_0)=g(\bold e_2,\check\bold F(\bold e_1))
-\lambda_0\,g(\bold e_2,\bold e_1)=\\
\vspace{1ex}
=g(\check\bold F(\bold e_2),\bold e_1)-\lambda_0\,g(\bold e_2,\bold e_1)
=\lambda_0\,g(\bold e_2,\bold e_1)-\lambda_0\,g(\bold e_2,\bold e_1)=0.
\endgathered
\mytag{5.16}
$$
According to the formulas \mythetag{5.15} and \mythetag{5.16} the matrix 
of the metric $\bold g$ in the basis $\bold e_0,\,\bold e_1,\,\bold e_2$ 
takes the following form:
$$
\hskip -2em
g_{ij}=\Vmatrix\format\l\ \ &\r &\ \ \r\\
0 & g_{01} & 0\\ 
\vspace{1ex} g_{01} &g_{11} & g_{12}\\ 
\vspace{1ex} 0 & g_{12} & g_{22}\endVmatrix.
\mytag{5.17}
$$
Note that $g_{22}\neq 0$ and $g_{01}\neq 0$ in \mythetag{5.17}. Indeed, 
otherwise the metric $\bold g$ would be degenerate. Note also that the
choice of the vectors $\bold e_0$, $\bold e_1$, and $\bold e_2$ is not
unique. We can perform the following transformations of these vectors:
$$
\hskip -2em
\aligned
&\bold e_0\to\alpha\,\bold e_0,\\
\vspace{1ex}
&\bold e_1\to\beta\,\bold e_0+\alpha\,\bold e_1+\gamma\,\bold e_2,\\
\vspace{1ex}
&\bold e_2\to\theta\,\bold e_0+\delta\,\bold e_2.
\endaligned
\mytag{5.18}
$$
The basis transformations of the form \mythetag{5.18} preserve the
relationships \mythetag{5.12}, \mythetag{5.13}, and \mythetag{5.14}.
Hence, the zero components of the matrix \mythetag{5.17} remain zero 
under these transformations. For the beginning we choose the special
transformation of the form \mythetag{5.18} by setting
$$
\hskip -2em
\alpha=1,\qquad\delta=1,\qquad\gamma=0,
\qquad\theta=-\frac{g_{12}}{g_{01}},
\qquad\beta=-\frac{g_{11}}{2\,g_{01}}
\mytag{5.19}
$$
in \mythetag{5.18}.
Applying the basis transformations \mythetag{5.18} with the parameters 
\mythetag{5.19}, we find that the components $g_{11}$ and $g_{12}$ in 
\mythetag{5.17} turn to zero. \pagebreak As a result the matrix of the 
metric $\bold g$ takes the following form:
$$
\hskip -2em
g_{ij}=
\Vmatrix\format\c\ \ &\c &\ \ \c\\
0 & g_{01} & 0\\ 
\vspace{1ex} g_{01} & 0 & 0\\ 
\vspace{1ex} 0 & 0 & g_{22}\endVmatrix.
\mytag{5.20}
$$
As we already noted, $g_{22}\neq 0$. If $g_{22}>0$, then due to 
\mythetag{5.20} the signature of the metric $\bold g$ would be
$(+,+,-)$. But, actually, the signature of $\bold g$ is $(+,-,-)$.
For this reason we conclude that $g_{22}<0$. The component 
$g_{01}$ is also nonzero. It can be either positive or negative.
The sign of the component $g_{01}$ in \mythetag{5.20} is a new 
invariant of a pair of forms in the third case. We denote it 
$\sigma_3$:
$$
\hskip -2em
\sigma_3=\sign(g_{01}).
\mytag{5.21}
$$
The invariant \mythetag{5.21} is analogous to other integer invariants
$\sigma_0$, $\sigma_1$, and $\sigma_2$ consi\-dered above. Depending on 
the value of this invariant we specify the subcase two and the
subcase three within the third case.
\par
     The {\bf subcase two} within the third case of our classification 
scheme is fixed by the following values of the invariants $D_3$, $D_2$, 
$\sigma_2$ and $\sigma_3$:
$$
\xalignat 4
&\hskip -2em
D_3=0,
&&D_2=0,
&&\sigma_2=1,
&&\sigma_3=1.
\quad
\mytag{5.22}
\endxalignat
$$
If the conditions \mythetag{5.22} hold, the matrix of the metric 
$\bold g$ can be brought to
$$
\hskip -2em
g_{ij}=
\Vmatrix\format\c\ \ \ &\c &\ \ \!\c\\
0 & b^2 & \,0\\ 
\vspace{1ex} b^2 & 0 & \,0\\ 
\vspace{1ex} 0 & 0 & -c^2\endVmatrix.
\mytag{5.23}
$$
The matrix \mythetag{5.23} is not an ultimate canonical presentation 
for the metric $\bold g$. In order to simplify \mythetag{5.23} we
apply the following basis transformation:
$$
\hskip -2em
\aligned
&\bold e_0\to\frac{1}{b}\,\bold e_1,\\
\vspace{1ex}
&\bold e_1\to\frac{1}{b}\,\bold e_0,\\
\vspace{1ex}
&\bold e_2\to\frac{1}{c}\,\bold e_2.
\endaligned
\mytag{5.24}
$$
Upon applying the basis transformation \mythetag{5.24} we find that 
the matrices of the metrics $\bold g$ and $\check\bold g$ take the 
following canonical forms, where $a=\lambda_0$:
$$
\xalignat 2
&\hskip -2em
g_{ij}=\Vmatrix\format\l\ \ \ \ &\r &\ \ \r\\
0 & 1 & 0\\ 
\vspace{1ex} 1 &0 & 0\\ 
\vspace{1ex} 0 & 0 &-1\endVmatrix,
&&\check g_{ij}=\Vmatrix\format\l\ \ \ \ &\r &\ \ \r\\
1 & a & 0\\ 
\vspace{1ex} a &0 & 0\\ 
\vspace{1ex} 0 & 0 &-a\endVmatrix.
\mytag{5.25}
\endxalignat
$$
\mytheorem{5.2} If \,the conditions \mythetag{5.22} are fulfilled, then
the matrices of the metrics $\bold g$ and $\check\bold g$ take their
canonical forms \mythetag{5.25} in some basis.
\endproclaim
     The {\bf subcase three} differs from the subcase two by the value
of the invariant $\sigma_3$. Instead of the equalities \mythetag{5.22},
here we have
$$
\xalignat 4
&\hskip -2em
D_3=0,
&&D_2=0,
&&\sigma_2=1,
&&\sigma_3=-1.
\quad
\mytag{5.26}
\endxalignat
$$
Under the conditions \mythetag{5.26} the matrix \mythetag{5.20} specifies 
to
$$
\hskip -2em
g_{ij}=
\Vmatrix\format\c\ \ &\c &\ \ \!\c\\
\,0 &- b^2 & \,0\\ 
\vspace{1ex} -b^2 & \,0 & \,0\\ 
\vspace{1ex} \,0 & \,0 & -c^2\endVmatrix.
\mytag{5.27}
$$
In order to simplify \mythetag{5.27} we apply the following 
basis transformation:
$$
\hskip -2em
\aligned
&\bold e_0\to-\frac{1}{b}\,\bold e_0,\\
\vspace{1ex}
&\bold e_1\to\frac{1}{b}\,\bold e_1,\\
\vspace{1ex}
&\bold e_2\to\frac{1}{c}\,\bold e_2.
\endaligned
\mytag{5.28}
$$
Upon applying \mythetag{5.28} we find that the matrices of the 
metric $\bold g$ and $\check\bold g$ take the following canonical 
forms, where $a=\lambda_0$:
$$
\xalignat 2
&\hskip -2em
g_{ij}=\Vmatrix\format\l\ \ \ \ &\r &\ \ \r\\
0 & 1 & 0\\ 
\vspace{1ex} 1 &0 & 0\\ 
\vspace{1ex} 0 & 0 &-1\endVmatrix,
&&\check g_{ij}=\Vmatrix\format\l\ \ \,&\r &\,\ \ \r\\
0 & a & 0\\ 
\vspace{1ex} a &-1 & 0\\ 
\vspace{1ex} 0 & 0 &-a\endVmatrix.
\mytag{5.29}
\endxalignat
$$
\mytheorem{5.3} If \,the conditions \mythetag{5.26} are fulfilled, then
the matrices of the metrics $\bold g$ and $\check\bold g$ take their
canonical forms \mythetag{5.29} in some basis.
\endproclaim
    Now, again keeping the conditions \mythetag{5.1} unchanged, we 
increase by one the value of the integer invariant $\sigma_2$ in 
\mythetag{5.8}. As a result we get
$$
\hskip -2em
\sigma_2=2.
\mytag{5.30}
$$
Combining \mythetag{5.1} with \mythetag{5.30}, we write 
$$
\xalignat 3
&\hskip -2em
D_3=0,
&&D_2=0,
&&\sigma_2=2.
\quad
\mytag{5.31}
\endxalignat
$$
The {\bf subcase four} within the third case of our classification 
scheme is specified by the conditions \mythetag{5.31}. If the equalities 
\mythetag{5.31} hold, we have the following subspaces:
$$
\hskip -2em
\aligned
&W=\Img((\check\bold F-\lambda_0\,I)^2)=\Ker(\check\bold F-\lambda_0\,I),\\
\vspace{1ex}
&U=\Ker((\check\bold F-\lambda_0\,I)^2)=\Img(\check\bold F-\lambda_0\,I).
\endaligned
\mytag{5.32}
$$
The dimensions of these subspaces \mythetag{5.32} are also given by the  
formulas \mythetag{5.10}. Let's choose some vector $\bold e_2$ in $V$ 
such that $\bold e_2\notin U$. It is clear that $\bold e_2\neq 0$. Then 
we define the vector $\bold e_1$ by means of the formula
$$
\hskip -2em
\bold e_1=\check\bold F(\bold e_2)-\lambda_0\,\bold e_2.
\mytag{5.33}
$$
Due to \mythetag{5.32} from \mythetag{5.33} we derive that $\bold e_1\in U$,
$\bold e_1\neq 0$, and $\bold e_1\notin W$. Therefore we can define the
third vector $\bold e_0$ by means of the formula
$$
\hskip -2em
\bold e_0=\check\bold F(\bold e_1)-\lambda_0\,\bold e_1.
\mytag{5.34}
$$
The vector $\bold e_0$ belongs to the subspace $W$. Since $\bold e_1
\notin W$, the vector $\bold e_0$ is nonzero. From $\bold e_0\in W$ we 
derive that $\bold e_0$ is an eigenvector of the associated operator 
$\check\bold F$ corresponding to the eigenvalue $\lambda_0$ in
\mythetag{5.2}:
$$
\hskip -2em
\check\bold F(\bold e_0)=\lambda_0\,\bold e_0.
\mytag{5.35}
$$
The vectors $\bold e_0$, $\bold e_1$, and $\bold e_2$ are linearly
independent. They form a {\bf Jordan normal basis} for the operator
$\check\bold F$. Let's study the components of the metric $\bold g$
in this basis. Applying the formulas \mythetag{5.34}, \mythetag{5.35}, 
and \mythetag{2.2}, we derive
$$
\gathered
g_{00}=g(\bold e_0,\bold e_0)=g(\bold e_0,\check\bold F(\bold e_1))
-\lambda_0\,g(\bold e_0,\bold e_1)=\\
\vspace{1ex}
=g(\check\bold F(\bold e_0),\bold e_1)-\lambda_0\,g(\bold e_0,\bold e_1)
=\lambda_0\,g(\bold e_0,\bold e_1)-\lambda_0\,g(\bold e_0,\bold e_1)=0.
\endgathered
\mytag{5.36}
$$
Similarly, applying the formulas \mythetag{5.33}, \mythetag{5.35}, 
and \mythetag{2.2}, we derive
$$
\gathered
g_{10}=g_{01}=g(\bold e_0,\bold e_1)=g(\bold e_0,\check\bold F(\bold e_2))
-\lambda_0\,g(\bold e_0,\bold e_2)=\\
\vspace{1ex}
=g(\check\bold F(\bold e_0),\bold e_2)-\lambda_0\,g(\bold e_0,\bold e_2)
=\lambda_0\,g(\bold e_0,\bold e_2)-\lambda_0\,g(\bold e_0,\bold e_2)=0.
\endgathered
\mytag{5.37}
$$
And finally, applying the formulas \mythetag{5.33}, \mythetag{5.34}, 
and \mythetag{2.2}, we derive
$$
\gathered
g_{11}=g(\bold e_1,\bold e_1)=g(\bold e_1,\check\bold F(\bold e_2))
-\lambda_0\,g(\bold e_1,\bold e_2)=\\
\vspace{1ex}
=g(\check\bold F(\bold e_1),\bold e_2)-\lambda_0\,g(\bold e_1,\bold e_2)
=\lambda_0\,g(\bold e_1,\bold e_2)\,+\\
\vspace{1ex}
+\,g(\bold e_0,\bold e_2)-\lambda_0\,g(\bold e_0,\bold e_2)
=g(\bold e_0,\bold e_2)=g_{02}=g_{20}.
\endgathered
\mytag{5.38}
$$
Due to the formulas \mythetag{5.36}, \mythetag{5.37}, and \mythetag{5.38}, 
the matrix of the metric $\bold g$ in the basis $\bold e_0,\,\bold e_1,\,
\bold e_2$ takes the following form: 
$$
\hskip -2em
g_{ij}=
\Vmatrix\format\l\ \ &\c &\ \ \ \r\\
0 & 0 & g_{11}\\ 
\vspace{1ex} 0 & g_{11} & g_{12}\\ 
\vspace{1ex} g_{11} & g_{12} & g_{22}\endVmatrix.
\mytag{5.39}
$$
The component $g_{11}$ in the matrix \mythetag{5.39} is nonzero since 
otherwise the metric $\bold g$ would be degenerate. \pagebreak Note 
that the Jordan normal basis $\bold e_0,\,\bold e_1,\,\bold e_2$ of the
associated operator $\check\bold F$ is not unique. Applying the 
transformation
$$
\hskip -2em
\aligned
&\bold e_0\to\alpha\,\bold e_0,\\
\vspace{1ex}
&\bold e_1\to\beta\,\bold e_0+\alpha\,\bold e_1,
\vspace{1ex}
&\bold e_2\to\gamma\,\bold e_0+\beta\,\bold e_1+\alpha\,\bold e_2,
\endaligned
\mytag{5.40}
$$
we get another Jordan normal basis, where the metric $\bold g$ is 
presented by another matrix of the form \mythetag{5.39}. Since
$g_{11}\neq 0$ we can set
$$
\xalignat 3
&\hskip -2em
\alpha=1,
&&\beta=-\frac{g_{12}}{g_{11}},
&&\gamma=-\frac{g_{22}}{2\,g_{11}}+\frac{3\,(g_{12})^2}{8\,(g_{11})^2}.
\qquad
\mytag{5.41}
\endxalignat
$$
Applying \mythetag{5.41} to \mythetag{5.40}, we find that the metric 
$\bold g$ is given by the following skew-diagonal matrix in the new
Jordan normal basis of the operator $\check\bold F$:
$$
\hskip -2em
g_{ij}=
\Vmatrix\format\c\ \ &\c &\ \ \c\\
0 & 0 & g_{11}\\ 
\vspace{1ex} 0 & g_{11} & 0\\ 
\vspace{1ex} g_{11} & 0 & 0\endVmatrix.
\mytag{5.42}
$$
The component $g_{11}$ in \mythetag{5.42} is negative since otherwise,
if $g_{11}>0$, the signature of the metric $\bold g$ would be $(+,+,-)$, 
while actually it is $(+,-,-)$. Therefore, we can substitute 
$g_{11}=-b^2$ into \mythetag{5.42}. As a result we get the matrix 
$$
\hskip -2em
g_{ij}=
\Vmatrix\format\c\ \ &\c &\ \ \c\\
0 & 0 & -b^2\\ 
\vspace{1ex} 0 & -b^2 & 0\\ 
\vspace{1ex} -b^2 & 0 & 0\endVmatrix,
\mytag{5.43}
$$
where $b\neq 0$. The matrix \mythetag{5.43} is an intermediate
result. In order to get the ultimate result we consider the 
following basis transformation:
$$
\hskip -2em
\aligned
&\bold e_0\to -\frac{1}{b}\,\bold e_0,\\
\vspace{1ex}
&\bold e_1\to\frac{1}{b}\,\bold e_2,\\
\vspace{1ex}
&\bold e_2\to -\frac{1}{b}\,\bold e_1.
\endaligned
\mytag{5.44}
$$
Upon applying the basis transformation \mythetag{5.44} we find 
that the matrices of the metrics $\bold g$ and $\check\bold g$ 
take their canonical forms
$$
\xalignat 2
&\hskip -2em
g_{ij}=\Vmatrix\format\l\ \ \ \ &\r &\ \ \r\\
0 & 1 & 0\\ 
\vspace{1ex} 1 &0 & 0\\ 
\vspace{1ex} 0 & 0 &-1\endVmatrix,
&&\check g_{ij}=\Vmatrix\format\l\ \ \ \ &\r &\ \ \r\\
0 & a & 0\\ 
\vspace{1ex} a & 0 & 1\\ 
\vspace{1ex} 0 & 1 &-a\endVmatrix.
\mytag{5.45}
\endxalignat
$$
\mytheorem{5.4} If \,the conditions \mythetag{5.31} are fulfilled, 
then the matrices of the metrics $\bold g$ and $\check\bold g$ take 
their canonical forms \mythetag{5.45} in some basis.
\endproclaim
\head
6. Classification.
\endhead
    Thus, all possible cases are exhausted. The total number of them is ten.
All of these ten cases are given in the following three tables.
$$
\vcenter{\hsize=200pt
\offinterlineskip
\halign{\vrule#&\hfill\quad #\quad\hfill&\vrule#&\quad #\quad\hfill
&\vrule#\cr
\noalign{\hrule}
height 14pt depth 8pt & Condition &&\hfill Canonical presentation&\cr
\noalign{\hrule}
&$D_3>0$ 
&& $\vcenter{\vskip -1ex
$$\aligned &g_{ij}=\Vmatrix\format\l\ \ &\r &\ \ \r\\
1 & 0 & 0\\ \vspace{1ex} 0 &-1 & 0\\ \vspace{1ex} 0 & 0 &-1\endVmatrix,
\quad\check g_{ij}=\Vmatrix\format\l\ \ \ &\r &\ \ \ \r\\
a & 0 & 0\\ \vspace{1ex} 0 &b & 0\\ \vspace{1ex} 0 & 0 &c\endVmatrix,\\
\vspace{2ex}
&\text{where \ }a\neq -b\text{, \ }a\neq -c\text{, \ and \ }b\neq c
\endaligned$$
\vskip -1.5ex}$
&\cr
\noalign{\hrule}
&$\vcenter{\hsize=1.5cm
$$
\gathered
D_3=0,\\
D_2>0,\\
\sigma_0=1
\endgathered
$$}$
&& $\vcenter{\vskip -1ex
$$\aligned &g_{ij}=\Vmatrix\format\l\ \ &\r &\ \ \r\\
1 & 0 & 0\\ \vspace{1ex} 0 &-1 & 0\\ \vspace{1ex} 0 & 0 &-1\endVmatrix,
\quad\check g_{ij}=\Vmatrix\format\l\ \ \ &\r &\ \ \ \r\\
a & 0 & 0\\ \vspace{1ex} 0 &b & 0\\ \vspace{1ex} 0 & 0 &b\endVmatrix,\\
\vspace{2ex}
&\text{where \ }a\neq -b
\endaligned$$
\vskip -1.5ex}$
&\cr
\noalign{\hrule}}}
\mytag{6.1}
$$

$$
\vcenter{\hsize=200pt
\offinterlineskip
\halign{\vrule#&\hfill\quad #\quad\hfill&\vrule#&\quad #\quad\hfill
&\vrule#\cr
\noalign{\hrule}
height 14pt depth 8pt & Condition &&\hfill Canonical presentation&\cr
\noalign{\hrule}
&$D_3<0$ 
&& $\vcenter{\vskip -1ex
$$\aligned &g_{ij}=\Vmatrix\format\l\ \ &\r &\ \ \r\\
1 & 0 & 0\\ \vspace{1ex} 0 &-1 & 0\\ \vspace{1ex} 0 & 0 &-1\endVmatrix,
\quad\check g_{ij}=\Vmatrix\format\l\ \ &\r &\,\ \ \ \r\\
a & b & 0\\ 
\vspace{1ex} b &-a & 0\\ 
\vspace{1ex} 0 & 0 &c\endVmatrix,\\
\vspace{2ex}
&\text{where \ }b\neq 0
\endaligned$$
\vskip -1.5ex}$
&\cr
\noalign{\hrule}
&$\vcenter{\hsize=1.5cm
$$
\gathered
D_3=0,\\
D_2\neq 0,\\
\sigma_0=-1,\\
\sigma_1=0
\endgathered
$$}$
&& $\vcenter{\vskip -1ex
$$\aligned &g_{ij}=\Vmatrix\format\l\ \ &\r &\ \ \r\\
1 & 0 & 0\\ \vspace{1ex} 0 &-1 & 0\\ \vspace{1ex} 0 & 0 &-1\endVmatrix,
\quad
\check g_{ij}=\Vmatrix\format\l\ \ &\r &\,\ \ \ \r\\
a & 0 & 0\\ 
\vspace{1ex} 0 &-a & 0\\ 
\vspace{1ex} 0 & 0 &c\endVmatrix,\\
\vspace{2ex}
&\text{where \ }a\neq -c
\endaligned$$
\vskip -1.5ex}$
&\cr
\noalign{\hrule}
&$\vcenter{\hsize=1.5cm
$$
\gathered
D_3=0,\\
D_2\neq 0,\\
\sigma_0=-1,\\
\sigma_1=1
\endgathered
$$}$
&& $\vcenter{\vskip -1ex
$$\aligned &g_{ij}=\Vmatrix\format\l\ \ &\r &\ \ \r\\
1 & 0 & 0\\ \vspace{1ex} 0 &-1 & 0\\ \vspace{1ex} 0 & 0 &-1\endVmatrix,
\quad
\check g_{ij}=\Vmatrix\format\l\ \ \ \ &\r &\,\ \ \ \r\\
1 & a & 0\\ 
\vspace{1ex} a &0 & 0\\ 
\vspace{1ex} 0 & 0 &c\endVmatrix,\\
\vspace{2ex}
&\text{where \ }a\neq -c
\endaligned$$
\vskip -1.5ex}$
&\cr
\noalign{\hrule}
&$\vcenter{\hsize=1.5cm
$$
\gathered
D_3=0,\\
D_2\neq 0,\\
\sigma_0=-1,\\
\sigma_1=-1
\endgathered
$$}$
&& $\vcenter{\vskip -1ex
$$\aligned &g_{ij}=\Vmatrix\format\l\ \ &\r &\ \ \r\\
1 & 0 & 0\\ \vspace{1ex} 0 &-1 & 0\\ \vspace{1ex} 0 & 0 &-1\endVmatrix,
\quad
\check g_{ij}=\Vmatrix\format\l\ \ \,&\r &\,\ \ \ \r\\
0 & a & 0\\ 
\vspace{1ex} a &-1 & 0\\ 
\vspace{1ex} 0 & 0 &c\endVmatrix,\\
\vspace{2ex}
&\text{where \ }a\neq -c
\endaligned$$
\vskip -1.5ex}$
&\cr
\noalign{\hrule}}}
\mytag{6.2}
$$
$$
\vcenter{\hsize=200pt
\offinterlineskip
\halign{\vrule#&\hfill\quad #\quad\hfill&\vrule#&\quad #\quad\hfill
&\vrule#\cr
\noalign{\hrule}
height 14pt depth 8pt & Condition &&\hfill Canonical presentation&\cr
\noalign{\hrule}
&$\vcenter{\hsize=1.5cm
$$
\gathered
D_3=0,\\
D_2=0,\\
\sigma_2=0
\endgathered
$$}$
&& $\vcenter{\vskip -1ex
$$g_{ij}=\Vmatrix\format\l\ \ &\r &\ \ \r\\
1 & 0 & 0\\ \vspace{1ex} 0 &-1 & 0\\ \vspace{1ex} 0 & 0 &-1\endVmatrix,
\quad\check g_{ij}=\Vmatrix\format\l\ \ &\r &\ \ \r\\
a & 0 & 0\\ 
\vspace{1ex} 0 &-a & 0\\ 
\vspace{1ex} 0 & 0 &-a\endVmatrix
$$
\vskip -1.5ex}$
&\cr
\noalign{\hrule}
&$\vcenter{\hsize=1.5cm
$$
\gathered
D_3=0,\\
D_2=0,\\
\sigma_2=1,\\
\sigma_3=1
\endgathered
$$}$
&& $\vcenter{\vskip -1ex
$$g_{ij}=\Vmatrix\format\l\ \ \ \ &\r &\ \ \r\\
0 & 1 & 0\\ 
\vspace{1ex} 1 &0 & 0\\ 
\vspace{1ex} 0 & 0 &-1\endVmatrix,
\quad\check g_{ij}=\Vmatrix\format\l\ \ \ \ &\r &\ \ \r\\
1 & a & 0\\ 
\vspace{1ex} a &0 & 0\\ 
\vspace{1ex} 0 & 0 &-a\endVmatrix
$$
\vskip -1.5ex}$
&\cr
\noalign{\hrule}
&$\vcenter{\hsize=1.5cm
$$
\gathered
D_3=0,\\
D_2=0,\\
\sigma_2=1,\\
\sigma_3=-1
\endgathered
$$}$
&& $\vcenter{\vskip -1ex
$$g_{ij}=\Vmatrix\format\l\ \ \ \ &\r &\ \ \r\\
0 & 1 & 0\\ 
\vspace{1ex} 1 &0 & 0\\ 
\vspace{1ex} 0 & 0 &-1\endVmatrix,
\quad
\check g_{ij}=\Vmatrix\format\l\ \ \,&\r &\,\ \ \r\\
0 & a & 0\\ 
\vspace{1ex} a &-1 & 0\\ 
\vspace{1ex} 0 & 0 &-a\endVmatrix
$$
\vskip -1.5ex}$
&\cr
\noalign{\hrule}
&$\vcenter{\hsize=1.5cm
$$
\gathered
D_3=0,\\
D_2=0,\\
\sigma_2=2
\endgathered
$$}$
&& $\vcenter{\vskip -1ex
$$
g_{ij}=\Vmatrix\format\l\ \ \ \ &\r &\ \ \r\\
0 & 1 & 0\\ 
\vspace{1ex} 1 &0 & 0\\ 
\vspace{1ex} 0 & 0 &-1\endVmatrix,
\quad
\check g_{ij}=\Vmatrix\format\l\ \ \ \ &\r &\ \ \r\\
0 & a & 0\\ 
\vspace{1ex} a & 0 & 1\\ 
\vspace{1ex} 0 & 1 &-a\endVmatrix
$$
\vskip -1.5ex}$
&\cr
\noalign{\hrule}}}
\mytag{6.3}
$$
\vskip 1ex
Apart from the discriminants $D_3$ and $D_2$, in the first columns 
of the tables \mythetag{6.1}, \mythetag{6.2}, and  \mythetag{6.3}
we see some special invariants $\sigma_0$, $\sigma_1$, $\sigma_2$, 
and $\sigma_3$. Currently I don't know any explicit formulas for
expressing these invariants trough the components of the metrics
$\bold g$ and $\check\bold g$ in an arbitrary basis. However, I 
expect that such formulas can be derived.

\enddocument
\end